\RequirePackage{ifpdf}
\ifpdf 
\documentclass[pdftex]{sigma}
\else
\documentclass{sigma}
\fi

\def \C{\mathbb C}
\def \Z {\mathbb Z}
\def \N {\mathbb N}

\def \e {{\epsilon}}
\def \Ci {{C^\infty}}
\def \Cl {{C\ell}}

\numberwithin{equation}{section}

\begin{document}

\allowdisplaybreaks

\renewcommand{\thefootnote}{$\star$}

\renewcommand{\PaperNumber}{077}

\FirstPageHeading

\ShortArticleName{A Canonical Trace Associated with Certain Spectral Triples}

\ArticleName{A Canonical Trace Associated\\ with Certain Spectral Triples\footnote{This paper is a
contribution to the Special Issue ``Noncommutative Spaces and Fields''. The
full collection is available at
\href{http://www.emis.de/journals/SIGMA/noncommutative.html}{http://www.emis.de/journals/SIGMA/noncommutative.html}}}

\Author{Sylvie PAYCHA}

\AuthorNameForHeading{S. Paycha}

\Address{Laboratoire de Math\'ematiques, 63177 Aubi\`ere Cedex, France }
\Email{\href{mailto:sylvie.paycha@math.univ-bpclermont.fr}{sylvie.paycha@math.univ-bpclermont.fr}}

\URLaddress{\url{http://math.univ-bpclermont.fr/~paycha/}}

\ArticleDates{Received March 11, 2010, in f\/inal form September 13, 2010;  Published online September 29, 2010}

\Abstract{In the  abstract pseudodif\/ferential setup of
Connes and Moscovici, we prove a general
formula for the discrepancies of zeta-regularised traces  associated with
certain spectral triples,  and we introduce a
canonical trace on operators, whose order lies outside (minus) the dimension
spectrum of the spectral triple.}


\Classification{58B34; 58J42; 47B47}

\renewcommand{\thefootnote}{\arabic{footnote}}
\setcounter{footnote}{0}

\section{Introduction}

Connes and Moscovici's setup for abstract pseudodif\/ferential calculus \cite{CoMo} (see also \cite{Hig})
associated with a certain type of spectral triple $\left({\cal A}, {\cal H}, D\right)$ provides
a framework, in which the ordinary trace on trace-class operators can be
extended to a linear form on all abstract pseudodif\/ferential
operators using zeta regualarisation type methods. Here,  ${\cal A}$ is an
involutive algebra represented in a complex Hilbert space ${\cal H}$ and $D$ a
self-adjoint operator in ${\cal H}$.  This linear  extension, which depends on the operator $D$  used as a regulator in the zeta regularisation
procedure,  does not vanish on commutators as it can be seen from results of~\cite{CoMo}.

These constructions mimick the classical pseudodif\/ferential calculus setup on
a closed manifold, where  similar linear extensions
of the ordinary trace are built using the same zeta regularisation type procedure. On non-integer order operators, these linear extensions are independent of
the regulator  used in the zeta regularisation procedure. They def\/ine  the canonical
trace popu\-larised by Kontsevich and Vishik~\cite{KV},  which  vanishes
on non-integer order commutators of classical pseudodif\/ferential operators. The
canonical trace is unique, in so far as any linear form on  non-integer order
classical pseudodif\/ferential operators  which  vanishes
on non-integer order commutators, is proportional to the canonical trace \cite{MSS}.

In the abstract pseudodif\/ferential calculus framework, and  under a mild
reinforcement (see assumption (H) in Section~\ref{section4}  and assumptions (Dk) and (T) in Section~\ref{section5})  of the usual assumptions on a regular
spectral triple with discrete dimension spectrum Sd (used  to def\/ine zeta regularised type linear extensions of
the trace), we show that  a similar canonical linear form on operators whose
order lies outside the discrete set $-{\rm Sd}$  can be
def\/ined, and that this canonical linear form
vanishes on commutators of operators, whose order lies outside the discrete set
$-{\rm Sd}$.

 The assumptions (H), (Dk) and (T) we put on the spectral triple generalise
the usual assumptions one puts on spectral triples to ensure the existence and a
reasonable pole structure of
meromorphic extensions of traces of holomorphic families of  operators of the
type $a \vert
D\vert^{-z}$ with $a$ in ${\cal A}$. Our strengthened assumptions ensure
the existence and a reasonable pole structure of meromorphic extensions of  traces of  more general holomorphic families
$a(z)   \vert D\vert^{-z}$ where $a(z)$ is a holomorphic family in ${\cal
  A}$.\footnote{\label{footnote:holfam}A family $\{f(z)\}_{z\in \Omega}$  in a topological vector space $E$,
 parametrised by a complex domain $\Omega$,  is {\it  holomorphic} at
$z_0\in \Omega$ if
the corresponding function $f:\Omega\to E$  is   uniformly
complex-dif\/ferentiable on compact subsets in a~neighborhood of $z_0$.  When $f$
takes its values in a Banach space $E$,    the existence of a complex
derivative  implies (via  a  Cauchy formula) that a holomorphic
function is actually inf\/initely dif\/ferentiable and admits a~Taylor expansion which converges uniformly  on a  compact disk centered at
  $z_0$ (see e.g.~\cite[Theorem~8.1.7]{Hig}   and \cite[Section~XV.5.1]{DS}).
   Holomorphy of a Banach space valued function is therefore equivalent to
   its analyticity.}  Broadly speaking, strengthening the assumptions on the spectral triple
amounts to embedding the spectral triple $\left({\cal A}, {\cal
    H}, D\right)$ into a holomorphic family $\left({\cal A}(z), {\cal
    H}, D\right)$ of spectral triples such that ${\cal A}(0)={\cal A}$.  Under these strengthened  assumptions,   not only can one  tackle
traces of holomorphic families of the type  $A
\vert D\vert^{-z}$ where $A$ is an abstract pseudodif\/ferential operator, but
one can also deal with
more general holomorphic families (see Def\/inition~\ref{defn:holfamily} in
Section~\ref{section4})
\[
A(z)\simeq\sum_{j=0}^\infty
b_j(z) \vert D\vert^{\alpha(z)-j},
\]
 where $b_j(z)$ is a holomorphic family in
the algebra ${\cal B}$  generated by  elements  of the type $\delta^n(a(z))$ or the type $ \delta^n\left([D,
  a(z)]\right)$ with $  n$ varying in $ \Z_{\geq 0}$ and  $a(z)$ any holomorphic family in ${\cal
  A} $.

 We derive a  general formula (see Theorem~\ref{thm:ResAz})\footnote{We have set    $ [[a,b]]= [a, b]\cap \Z$ for two integers $a<b$, ${\rm Res}^j_0$ stands for the $j$-th order
  residue at zero and ${\rm Res}_0^0$ stands for the f\/inite part at $z=0$.}
\[
{\rm Res}^{j+1}_{0} {\rm Tr}(A(z))\vert^{\rm mer}= \sum_{n=0}^{k-j}
\frac{\tau_{j+n}^{\vert D\vert}\left(\left(\partial_z^n (A(z)
  \vert D\vert^{q  z})\right)_{\vert_{z=0}}\right)}{q^{j+n+1}  n!}\qquad \forall\,
  j\in[[-1,k]],
\]
which expresses the complex residues
at zero of the meromorphic extensions ${\rm Tr}(A(z))\vert^{\rm
  mer}$ of the trace Tr of a holomorphic
family $A(z)$  in terms of residue type linear forms (see formula
(\ref{eq:taukD}))
\[ \tau^{\vert D\vert}_j(A):={\rm Res}^{j+1}_{0} {\rm
    Tr} \left( A  \vert D\vert^{-z}\right)\vert^{\rm mer}
    \]
introduced by Connes and
Moscovici in~\cite{CoMo}. Here  $ {\rm
    Tr} \left( A\, \vert D\vert^{-z}\right)\vert^{\rm mer}$
stands for  the meromorphic  extension of the holomorphic function
    ${\rm Tr} \left( A\, \vert D\vert^{-z}\right)$   def\/ined      on an
    adequately chosen  half-plane.

     Specialising to appropriate holomorphic families
yields explicit formulae for the discrepancies of the linear extensions of the
ordinary trace.  With the help of a holomorphic family
$(\vert D\vert+P)^{-z}$ def\/ined in (\ref{eq:phiDpert}) where $P$ is a
zero order perturbation, we can change the weight $\vert D\vert$ to
 $\vert D\vert+P$ in the linear forms $\tau_j^{\vert D\vert}$
  in order to def\/ine the linear form $\tau_j^{\vert D\vert +P}$ (see formula~(\ref{eq:phiDpert})).
  The subsequent identity (see
(\ref{eq:taujD1D2}))
\[
\tau_j^{\vert D\vert +P}(A)-\tau_j^{\vert D\vert}(A)= \sum_{n=0}^{k-j}
\frac{\tau_{j+n}^{\vert D\vert}\left( A \left(\partial_z^n((\vert D\vert+P)^{-z} \vert
    D\vert^z)_{\vert_{z=0}}\right)\right)}{n!}\qquad \forall\,
  j\in[[-1,k]],
  \]
 measures the sensitivity  of
$\tau_j^{\vert D\vert}$ to zero order perturbations   of the weight
$\vert D\vert$ and the following identity
(see
(\ref{eq:taujDbracket}))\footnote{We have set $L(A):= [\log \vert D\vert, A]$.}
\[
\tau_j^{\vert D\vert}([A,B])=\sum_{n= 0}^{k-j}
 (-1)^{n+1}\frac{\tau_{j+n}^{\vert D\vert}\left(
    A\,L^n(B)\right)_{\vert_{z=0}}}{n!},\qquad \forall\,
  j\in[[-1,k]],
  \]  which was derived in
Proposition~II.1 in~\cite{CoMo} by other means,  measures the obstruction to
its vanishing on commutators.

On the grounds of these identities, under assumptions (H), (Dk) and (T) we show that the  linear form
$\tau_{-1}^{\vert D\vert}$  is invariant under zero order perturbations of
$\vert D\vert$ (see Proposition~\ref{proposition6}) and that it vanishes on commutators for operators
whose order lies outside minus the dimension spectrum (see Proposition~\ref{proposition7}). These remarkable features  allow us to
consider this linear form as a substitute in the abstract pseudodif\/ferential operator set up
for Kontsevish and Vishik's canonical trace.

 To conclude, at the cost of
embedding a spectral triple in a holomorphic family of spectral triples, we
have built  a  linear form on  operators whose order lies outside a
given discrete set,  which  vanishes on commutators, whose orders lie outside
that discrete set and which is   insensitive to zero order perturbations of the weight $\vert D\vert$
used to build
this linear form.

\section{Spectral triples and abstract dif\/ferential operators}\label{section2}

We brief\/ly recall the setup of abstract dif\/ferential calculus as
explicited in~\cite{Hig}. It encompasses essential features of ordinary
pseudodif\/ferential calculus on manifolds, as  illustrated on a typical
example throughout this paragraph.

 Starting from a Hilbert space $ {\cal H}$ with scalar product $\langle\cdot, \cdot\rangle$ and associated norm $\Vert \cdot\Vert$ together with   a
self-adjoint operator $D$ on ${\cal H}$, we build a
non-negative self-adjoint operator $\Delta= D^2$.  For simplicity,  {\it we assume
the operator $\Delta$ is
invertible and hence positive}; if this is not the case we replace it  with
$1+\Delta$.
\begin{example}\label{ex:spin} Let $M$ be a closed and  smooth spin Riemannian  manifold, let  $E$
 be  a Clif\/ford Hermitian bundle
  over $M$ and $\nabla^E$ a Clif\/ford connection on $E$. Let $\Ci(M,E)$ denote
  the space of smooth sections of $E$ and  its $L^2$-closure ${\cal H}=L^2(M, E)$ for
  the Hermitian product induced by the metric on $M$ combined with the Hermitian
  structure on the f\/ibres of $E$. To this data corresponds a generalised Dirac
  operator (see e.g.~\cite{BGV})
  \[
  D= \sum_{i=1}^n c(dx_i)
  \nabla_{\frac{d}{dx_i}}^E.
  \]
\end{example}

 The self-adjoint operator $\Delta$ and its powers $\Delta^k$ are def\/ined on
 dense  domains    $ {\rm
    Dom}(\Delta^k)={\rm
    Dom}(D^{2k})$ in ${\cal H}$ and we consider the set
    \[
    {\cal H}^\infty:= \cap_{k=1}^\infty {\rm
    Dom}(\Delta^k)=\cap_{k=1}^\infty {\rm
    Dom}(D^k).
    \]

\begin{example}
With the notation of Example~\ref{ex:spin},  the operators
$\Delta^k=D^{2k}$ acting on the space of smooth sections of $E$ are essentially self-adjoint and have  unique self-adjoint extensions (denoted by the same symbol $\Delta^k$) to the   $H^k$-Sobolev spaces
  ${\rm Dom}(\Delta^k)=H^{2k}(M, E)$, obtained as closures of the space
  $\Ci(M,E)$ of smooth sections of the bundle~$E$ for the Sobolev norm
\[
\Vert u\Vert_{k}:= \big(\Vert u\Vert^2+ \Vert \Delta^{\frac{k}{2}}  u\Vert^2\big)^{\frac{1}{2}}
 \]
 (where $\Vert\cdot\Vert$ stands for the norm
on ${\cal H}$), and whose intersection ${\cal H}^\infty$
  coincides with  $\Ci(M, E)$.
\end{example}

We now introduce spectral triples~\cite{Co}, which  are the building blocks for abstract pseudodif\/ferential operators~\cite{CoMo,Hig}.

\begin{definition} A triple   $({\cal A}, {\cal H}, D)$, with  ${\cal A}$  an
involutive algebra represented in a complex Hilbert space ${\cal H}$ and $D$ a
self-adjoint operator in ${\cal H}$, is a {\it spectral triple}, if it satisf\/ies the
following properties:{\samepage
\begin{enumerate}
\itemsep=-0.5pt
\item[1)] the operator $a  (1+\Delta)^{-1}$ is compact for every element $a$ in
  ${\cal A}$,
\item[2)] the domain ${\rm Dom}(D)$ is stable under left multiplication by
${\cal A}$,
\item[3)] the adjoint operator $\delta_D: a\mapsto [D, a]$ on ${\cal A}$  extends to a bounded operator on ${\cal H}$.
\end{enumerate}}
\end{definition}

 For any complex number $s$
we set $\Delta^s:= \int_0^\infty
 \lambda^s   dP_\lambda$, where $P_\lambda:= 1_{({-}\infty, \lambda]}(\Delta) $
 stands for the spectral projection
of  the self-adjoint operator $\Delta$ corresponding to the interval $({-}\infty,
\lambda]$. The space ${\cal H}^s$ is obtained as the completion of ${\cal H}^\infty$
 for the norm\footnote{Recall that $\Delta$ is assumed to be invertible; if this
 is not the case we use $1+\Delta$ instead.}
 \[
 \Vert u\Vert_s:=\Vert \Delta^{\frac{s}{2d}} u\Vert.
 \]
  If  ${\rm Re}(s)\geq 0$, then  ${\cal H}^s$ lies in $ {\cal H}$ and coincides with
   the domain $ {\rm
  Dom}\big(\Delta^{\frac{s}{2}}\big)$.   If  ${\rm Re}(s)< 0$, then
  $\Delta^{\frac{s}{2}}$ is bounded on ${\cal H}$, whereas we have
  $  {\cal H} \subsetneq {\cal H}^s$.

  In the sequel, we shall often use the operator $\vert D\vert:= \sqrt{\Delta}$
  and  the unbounded derivation ${\rm ad}_{\vert D\vert}:
  X\mapsto [\vert D\vert, X]$ of the algebra  ${\cal B}({\cal H})$ of bounded
  linear operators on ${\cal H}$. Its
     domain is given by the set of bounded operators  that map the domain of
    $\vert D\vert$  onto itself, and for  which the commutator extends
    to a bounded operator on ${\cal H}$.

\begin{definition}
The spectral triple is said to be {\it regular} if and only if ${\cal A}$ and
$[D, {\cal A}]$ belong to the domain of smoothness\footnote{These notations are borrowed from \cite{CoMo},
  see Appendix B.}
  \[
  {\rm Op}^0:=\cap_{n=1}^\infty  {\rm Dom}(\delta^n),
   \]
   of the derivation $\delta$ on the algebra ${\cal B}({\cal H})$.
 \end{definition}

\begin{example}  With the notation of Example~\ref{ex:spin}, $\left(\Ci(M, E),
L^2(M, E), D\right)$ is a regular spectral triple.
\end{example}

For any complex number $s$, let ${\rm Op}^s$ be the set
of operators def\/ined by:
\[
P\in {\rm Op}^s \Longleftrightarrow \vert D\vert^{-s} P\in {\rm Op}^0 .
\]

Let us also introduce the  algebra\footnote{Here $\langle\cdots \rangle$
means ``generated by''.}
 \begin{gather*}
 {\cal B}=
 \langle \delta^{n_i}(a_i),\  \delta^{n_j}\left([D, a_j]\right), \  n_i, n_j\in \Z_{\geq 0}, \  a_i, a_j\in {\cal
  A}\rangle.
  \end{gather*}
 For a regular spectral triple $({\cal A}, {\cal H}, D)$, we have ${\cal B}\subset {\rm Op}^0$
and  the operator $b  \vert D\vert^s$ lies in ${\rm Op}^s$ for any $b$
in ${\cal B}$.

{\it We henceforth assume that  the spectral triple is regular and that ${\cal H}^\infty$ is stable under left
multiplication by ${\cal A}$.}

Following Higson \cite{Hig} (see Def\/inition 4.28 for the integer order case), for a complex number~$a$ we consider the set
${\cal E}^a\left({\cal A},D\right)$ of
operators $A:{\cal H}^\infty\to {\cal H}^\infty$, called  {\it  basic pseudodif\/ferential operators of order
 $a$}, that have the following expansion:
\[
A\simeq \sum_{j=0}^\infty b_j \vert D\vert^{a-j}, \qquad b_j\in {\cal
  B},
\]
   meaning by this, that for any non-negative integer $N$
\[
R_N(A):=A-\sum_{j=0}^{N-1}  b_{j} \vert D\vert^{a-j} \in {\rm  Op}^{a -N}.
\]
Any f\/inite linear combination of basic pseudodif\/ferential operators of order
$a$  is called
an {\it  abstract pseudodif\/ferential operator} of  order $a $.
By Proposition~4.31 in~\cite{Hig}
and Appendix~B in~\cite{CoMo}, abstract pseudodif\/ferential operators of integer order form an
algebra f\/iltered by the order, which is stable under the adjoint operator
$\delta=[\vert D\vert, \cdot]$.

\begin{example} With the notation  of Example \ref{ex:spin}, the triple consisting of  ${\cal A}=\Ci(M)$,
  ${\cal H}= L^2(M, E)$ and $D$ a generalised Dirac operator acting on $\Ci(M,
  E)$ form a regular spectral triple. Abstract pseudodif\/ferential operators
 of integer order correspond to a subalgebra of the algebra $\Cl^\Z(M, E)$ of classical pseudodif\/ferential operators of integer
  order \cite{Ho, Sh,Ta,Tr}, acting on smooth sections of the vector bundle $E$ over $M$. Since
   $D$ is a priori not an abstract pseudo\-dif\/ferential operator, the subalgebra
  corresponding to all integer order abstract pseudodif\/ferential operators is a
  priori smaller than  $\Cl^\Z(M, E)$.
\end{example}

By Lemma
4.30 in \cite{Hig} (see also (11) in Part II of \cite{CoMo}), for any abstract
pseudodif\/ferential operator $A$ and any complex number $\alpha$  we have:
 \begin{gather} \label{eq:commDalpha} \vert D\vert^\alpha\,  A\simeq    A\,
 \vert D\vert^\alpha+  \sum_{k=1}^\infty c_{\alpha,k}
\delta^k(A) \, \vert D\vert^{\alpha -k},
\end{gather}
 where we have set $c_{ \alpha, k}=
 \frac{\alpha(\alpha-1)\cdots (\alpha-k+1)}{k!}$ for any positive integer $k$. Note that if $A$ has order
$a$,
 then  $\delta^k(A)\, \vert D\vert^{\alpha-k}$ has order $a+\alpha-k$,   whose
 real value  decreases as $k$ increases.
\begin{example}
This formula holds  for any classical pseudodif\/ferential
operator $A$ acting on smooth sections of a vector bundle $E\to M$ over a
closed manifold $M$  and $D$ a generalised Dirac operator as in Example~\ref{ex:spin}. In this case,   the
limit $[\vert D\vert^\alpha, A]= \lim\limits_{N\to \infty}  \sum\limits_{k=1}^N c_{\alpha,k}
\delta^k(A)   \vert D\vert^{\alpha -k}$  holds in the Fr\'echet topology  of
classical pseudodif\/ferential operators with constant order (see~\cite{KV}),
 here $a+\alpha$ where $a$ is the order of~$A$.
\end{example}

\begin{definition} A  spectral triple $\left({\cal A}, {\cal H}, D\right)$ is   {\it
    finitely summable}   if  for any element $ a$ in
${\cal A}$, the operator $a\, \vert D\vert^{-1}$ lies in  a  Schatten class ${\cal
  L}^p({\cal H})$ for some $1\leq p<\infty$. Let $n$ stand for the inf\/imum of
the values  $p$ for which this holds, called {\it the degree of summability}
of $D$.
\end{definition}

\begin{remark}  In particular, for any element $ a$ in
${\cal A}$, the operator $a
    \vert D\vert^{-z}$ is trace-class  on the half-plane ${\rm
    Re}(z)>n $. Holomorphicity of the map $z\mapsto {\rm Tr}\left(a
    \vert D\vert^{-z}\right)$ on the half-plane ${\rm
    Re}(z)>n $ then follows. Indeed,  on any half-plane  ${\rm    Re}(z)>n+\e$ for some positive $\e$, the operator   $a   \vert D\vert^{-z}$  can be written  as  the product of a f\/ixed trace-class operator $A:= a \vert D\vert^{-n-\e}$ and a holomorphic family of bounded operators $B(z):=\vert D\vert^{-z+n+\e}$  on that
half-plane.   Since  holomorhicity implies analyticity for  Banach spaces valued functions (see  footnote~\ref{footnote:holfam}), there are bounded operators~$B_n$, $n\in \Z_{\geq 0}$ such that
$B(z)=\sum\limits_{n=0}^\infty B_n z^n$ converges uniformly on any compact disk centered at a~point of the half-plane  ${\rm
    Re}(z)>n+\e$.  Since $A$ is trace-class, so are the operators $AB_n$ and it  follows that ${\rm Tr}(a   \vert D\vert^{-z})= \sum\limits_{n=0}^\infty{\rm Tr}( A  B_n)z^n$ is holomorphic on every half-plane  ${\rm
    Re}(z)>n+\e$ for  positive~$\e$
 and hence on the half-plane  ${\rm
    Re}(z)>n$.
\end{remark}

More generally, given a f\/initely summable regular spectral triple  $\left({\cal A}, {\cal H},
  D\right)$ with degree of summability $n$, for any $A $ in ${\cal E}^a({\cal A}, D)$, the map $z\mapsto {\rm Tr}(A \,\vert D\vert^{-z})\vert^{\rm
  mer}$ is holomorphic on some half-plane ${\rm Re}(z-a)>n$ in  $\C$
depending on the order~$a$ of~$A$.

\begin{definition}
A regular and f\/initely summable spectral triple  has   {\it  discrete dimension
spectrum}, if  there is a discrete subset $S\subset \C$ such that for any
operator $A\in {\cal E}^a({\cal A}, D)$ with order $a$, the  map
 $z\mapsto {\rm Tr}(A  \vert D\vert^{-z})$
extends to a meromorphic map  $z\mapsto
 {\rm Tr}(A  \vert D\vert^{-z})\vert^{\rm
  mer}$ on $\C$
with  poles   in the  set
$S-a$. Let ${\rm Sd}$ denote the smallest such set $S$ called the {\it dimension
spectrum}.

 The dimension spectrum is {\it simple} if all the poles are
simple. It has f\/inite multiplicity  $k\in \N$ if the poles are at most of order $k$.\end{definition}

\begin{example} With the notation of   Example~\ref{ex:spin}, let $A$ be a
  classical pseudodif\/ferential operator of real order $a$, acting on
  $\Ci(M, E)$. The operator $z\mapsto {\rm Tr}(A  \vert
   D\vert^{-z})$  is holomorphic on  the half-plane ${\rm Re}(z)>a+n$, where
   $n$ is the dimension of the underlying manifold $M$, and has a~meromorphic
   extension to the whole complex plane with simple poles  in $]{-}\infty,
   a+n]\cap \Z$ (see~\cite{Se}).   This meromorphic extension can be written
   in terms of the canonical trace popularised by Kontsevich and Vishik \cite{KV} (see also~\cite{L})
   \[
   z\mapsto {\rm TR}(A  \vert    D\vert^{-z}).
   \]   The regular spectral triple $\left({\cal A}, {\cal H}, D\right)$
 arising from  Example~\ref{ex:spin}
 therefore has a  discrete simple dimension spectrum given by
  ${\rm Sd}= ]{-}\infty, n]\cap \Z$, which  is stable
under translations by negative integers. Note that the meromorphic extension is
holomorphic at zero if the order $a$ is non-integer.
\end{example}

\section{Logarithms}\label{section3}

Given a regular spectral triple\footnote{Recall that $\Delta$ is assumed to be
  invertible; otherwise we replace $\Delta$ with $\Delta+1$.}
$({\cal A}, {\cal H}, D)$,    the logarithm being a continuous
  function on the spectrum of the
operator
$\vert D\vert=\sqrt{\Delta}$,  one can
def\/ine the
unbounded self-adjoint operator $\log \vert D\vert$ by Borel functional
 calculus. It can  also be viewed as
the derivative at zero of the complex power $\vert D\vert^z$.  For any positive $\e$, the map  $z\mapsto \vert D\vert^{z-\e}\in {\rm Op}^{z-\e}$
 def\/ines a holomorphic  function on the half plane ${\rm Re}(z)<\e$ with
values  in ${\cal B}\left({\cal H}\right)$ and we have
$\log \vert D\vert =\vert D\vert^{\e}
 \left(\partial_z \left( \vert D\vert^{z-\e}\right)\right)_{\vert_{z=0}}
 $.
For any positive $\e$, the operator $ \log \vert D\vert   \vert D\vert^{-\e}=
\vert D\vert^{-\e}  \log \vert D\vert$ lies in   ${\rm Op}^0$, so that   $\log
\vert D\vert$ lies  in ${\rm Op}^\e$ for any positive~$\e$.

 For any complex number $a$ and any positive integer $L$,  we now introduce
 the set
${\cal E}^{a, L}({\cal A}, D) $   of
operators $A:{\cal H}^\infty\to {\cal H}^\infty$, which have the following expansion:
\[
A\simeq \sum_{l=0}^L\sum_{j=0}^\infty b_{j,l} \vert D\vert^{\alpha-j}  \log^l \vert D\vert, \qquad b_{j, l}\in {\cal
  B},
  \]
   meaning by this that for any non-negative integer $N$
\[
R_N(A):=A-\sum_{l=0}^L\sum_{j=0}^{N-1}  b_{j,l} \vert D\vert^{\alpha-j}  \log^l \vert D\vert  \in {\rm  Op}^{\alpha -N+\e}, \qquad\forall \,  \e>0.
\]
By convention we set  ${\cal E}^{\alpha, 0}({\cal A}, D)=
 {\cal E}^\alpha({\cal A}, D)$.

 \begin{remark} In the setup of pseudodif\/ferential
 operators acting on smooth sections of a vector bundle over a closed manifold
 $M$, for an operator $D$ as in Example \ref{ex:spin}, the set ${\cal E}^{a,
 L}({\cal A}, D)$ is reminiscent of the set $\Cl^{a,L}(M,E)$ of
 log-polyhomogeneous operators of order $a$ and logarithmic type $L$
 investigated in \cite{L}. Setting $L=0$ yields back classical
 pseudodif\/ferential operators.
 \end{remark}
We call  a linear operator $A: {\cal
   H}^\infty \to {\cal H}^\infty$ in ${\cal E}^{a, L}({\cal A}, D)$,  a {\it   basic    pseudodif\/ferential
   opera\-tor of order
 $a$ and   logarithmic  type $L$.} An abstract pseudodif\/ferential operator of
 order $a$ and logarithmic type $L$ is a f\/inite linear combination of basic
 abstract pseudodif\/ferential operators of order $a$ and logarithmic type $L$.

 Setting $a=z$ in~(\ref{eq:commDalpha}) and dif\/ferentiating with respect to $z$
 at zero yields
\begin{gather*}  \log \vert D\vert   A\simeq A   \log \vert D\vert +
 \sum_{k=1}^\infty c^\prime_{0,k}
\delta^k(A)   \vert D\vert^{-k},
\end{gather*}
for any abstract pseudodif\/ferential operator $A$,
with
$c^\prime_{0, k}:=\partial_z c_{ z, k}\vert_{z=0}=
\sum\limits_{j=1}^{k-1}\frac{(-1)^j}{j}$ for positive integers $k$. In particular,
this implies that the logarithmic type does not increase
 by the adjoint action with $ \log \vert D\vert$,
 so that for any complex number $a$ and any integer $k$,  we have
\begin{gather}\label{eq:commDalphalog}
A\in {\cal E}^{a,k}({\cal A}, D)\Longrightarrow [ \log \vert D\vert,    A ]
\in {\cal E}^{a-1,k}({\cal A}, D).
\end{gather}
\begin{remark}
  This  is reminiscent of   the  fact that in the classical setup,
  the  bracket of the logarithms of an elliptic operator
   (with appropriate conditions on the spectrum for the logarithm to
    be def\/ined) with a classical pseudodif\/ferential operator, is
    classical in spite of the fact that the   logarithm is not classical.
    \end{remark}

     The following proposition compares  the logarithms
   of two operators of the  same order.
\begin{proposition}\label{prop:differencelog}
For any $P\simeq \sum\limits_{j=0}^\infty b_j   \vert D\vert^{-j}
\in {\cal E}^0({\cal A}, D)$, the logarithm $\log
  (\vert D\vert+P)$ is well-defined and  lies in
${\cal E}^{0,1}({\cal A}, D)$. More precisely, the difference
 $\log   (\vert D\vert+P) -\log  \vert D\vert $ of the
  two logarithms lies  in ${\cal E}^0({\cal A}, D)$.
\end{proposition}

\begin{proof} By  (\ref{eq:commDalpha}) applied to $\alpha=-1$, the operator
$\vert D\vert^{-1}  P$ is an abstract pseudodif\/ferential operator of order $-1$
and hence a bounded operator on ${\cal H}$.   The operator $\log   (1+\vert
D\vert^{-1}  P)$
    def\/ined by  Borel functional calculus reads:
\[
\log  (1+\vert D\vert^{-1}  P)\simeq \sum_{i=0}^\infty
\frac{(-1)^{i-1}}{i!} (\vert D\vert^{-1}\, P)^i,
\] and hence lies in
   ${\cal E}^0({\cal A}, D)$ by (\ref{eq:commDalpha}).

The Campbell--Hausdorf\/f formula (see \cite{O} in the classical case) then yields:
\begin{gather*} \log   (\vert D\vert+P) -\log  \vert D\vert  \simeq
 \log (1+ \vert D\vert^{-1}  P) +\sum_{k=2}^{\infty } C^{(k)}(
  (\log \vert D\vert, 1+ \vert D\vert^{-1}  P) )
  ),
\end{gather*}where  $C^{(k)}( \cdot , \cdot)$ are  Lie monomials given by:
\[
C^{(k)}(A,B):=\sum\limits_{j=1}^\infty \frac{(-1)^{j+1}}{(j+1)} \sum
\frac{({\rm ad }_A)^{\alpha_1}({\rm ad }_B)^{\beta_1}\cdots
 ({\rm ad }_A)^{\alpha_j}({\rm ad }_B)^{\beta_j}B}{(1+ \sum_{l=1}^j
  \beta_l)\\  \alpha_1!\cdots\alpha_j! \beta_1!\cdots\beta_j!},
 \]
  with the inner sum running over $j-$ tuples of pairs $(\alpha_i,
  \beta_i)$ such that $\alpha_i\!+\!\beta_i\!>\!0$ and $\sum\limits_{i=1}^j
  \alpha_i\!+\! \beta_i\!=\! k$, and with the following notational convention:
   $ ({\rm ad
  }_A)^{\alpha_j}({\rm ad }_B)^{\beta_j-1} B=  ({\rm ad
  }_A)^{\alpha_j-1}  A$ if $\beta_j=0$,  in which case this vanishes if
  $\alpha_j>1$. Here ${\rm ad}_X:=[X, \cdot]$ denotes the adjoint
  action  by an opera\-tor~$X$. Property (\ref{eq:commDalphalog}) applied to $A=
  \log(1+ \vert D\vert^{-1}\,P)$,
   yields that
   \[
    C^{(1)}\left(\log \vert D \vert,
  \log(1+ \vert D\vert^{-1} P)\right)= \frac{1}{2}
\left[ \log \vert D \vert,  \log(1+ \vert D\vert^{-1} P)\right]
\] lies in
${\cal E}^0({\cal A}, D)$. By induction on~$k$, one shows that
$
C^{(k)}\left(\log \vert D \vert,
  \log(1+ \vert D\vert^{-1} P)\right)
$  lies in
${\cal E}^0({\cal A}, D)$ for any positive integer $k$,  so that   the
  dif\/ference $\log (\vert D\vert+P)-\log \vert D\vert $ also lies in ${\cal
  E}^0({\cal A}, D)$.

 Note that by (\ref{eq:commDalpha}) the adjoint operation
    ${\rm ad }_{\log(1+ \vert D\vert^{-1}\,B)} $  decreases the order by $1$
     unit and that   by~(\ref{eq:commDalphalog}) the same property holds
      for  the adjoint operation  ${\rm ad }_{\log \vert D\vert} $.  \end{proof}

\begin{remark}
Proposition \ref{prop:differencelog}  is reminiscent of   the  fact that in the classical setup,
  the dif\/ference of the logarithms of two elliptic operators
  (with appropriate conditions on the spectrum for their
   logarithms to be def\/ined) is classical in spite of the fact that each
    logarithm is not classical.\end{remark}

\section{Holomorphic families of operators}\label{section4}

 Let  $\left({\cal A}, {\cal H}, D\right)$ be a regular spectral triple  with
 discrete f\/inite dimension spectrum ${\rm Sd}$ and deg\-ree of summability $n$.
 Under the regularity assumption, one can equip ${\cal A}$ with a locally convex topology by means of semi-norms $a\mapsto \Vert \delta^k(a)\Vert$ and $a\mapsto \left \Vert \delta^k\left( [D, a]\right)\right\Vert$. The completion of~${\cal A}$ for this topology is a  Fr\'echet space~\cite{V} (and even a Fr\'echet $C^*$-algebra), which we denote by the same symbol ${\cal A}$, so that we shall henceforth assume that ${\cal A}$ is a Fr\'echet space.   Equivalence between
 holomorphicity and analyticity for   Banach space valued functions (see footnote~\ref{footnote:holfam})  actually extends to the case of Fr\'echet valued functions; the proof in the Banach case, which uses the Cauchy formula indeed extends to a Fr\'echet space in replacing the norm
  on the Banach space by each of the  semi-norms
  that def\/ine the topology on the Fr\'echet space.

   We call a ${\cal B}$-valued function $b(z)$
{\it holomorphic}, if it is a linear combination of products of  operators
$ \delta^{n_i}(a_i(z))$ and $\delta^{n_j}([\vert D\vert, a_j(z)])$
 for  some non-negative integers $n_i$ and $n_j$, where  $a_i(z)$
 and $a_j(z)$ are holomorphic families in ${\cal A}$.
 \begin{itemize}\itemsep=0pt
\item {\bf Assumption (H).}  {\it For any holomorphic family  $z\mapsto b(z)$ in ${\cal B}$
the  map $ z\mapsto {\rm Tr}(b(z) \vert D\vert^{-z})$
is holomorphic on the  half plane ${\rm Re}(z)>n$ and
 extends to a meromorphic map
 \begin{gather}
 \label{eq:TRz}
  z\mapsto {\rm Tr}(b(z) \vert D\vert^{-z})\vert^{\rm mer},
 \end{gather}
with  poles   in the same  set ${\rm Sd}$.}
\end{itemize}

\begin{example}
  With the notation Example \ref{ex:spin}, the spectral triple $\{\Ci(M), L^2(M, E), D\}$ sa\-tisf\/ies Assumption~(H). This follows from
  the fact that for any
  holomorphic family $z\mapsto f(z)$ in $\Ci(M)$, seen as a family of
  multiplication operators, the map $z\mapsto {\rm Tr}\big(f(z)   (1+
    D^2 )^{-z/2}\big)$ has a meromorphic extension  $z\mapsto {\rm TR}\big(f(z)   (1+
    D^2 )^{-z/2}\big)$ with simple poles in the set $ \Z\cap
]{-}\infty,n]$ and where, as before, ${\rm  TR}$ stands for the canonical
trace.
\end{example} This is actually a particular instance (with  $A(z)= f(z)   ( 1+
    D^2 )^{-z/2}$) of the following  more general result,
which can be found in~\cite{KV} (see also~\cite{Pa1} for a review). For any holomorphic family of classical pseudodif\/ferential
operators~$A(z)$ with non-constant af\/f\/ine order~$\alpha(z)$ acting on
$\Ci(M,E)$, the map  $z\mapsto {\rm TR}\left(A(z)\right)$ is meromorphic with
simple poles in the discrete set $\alpha^{-1}\left( \Z\cap   [{-}n, +\infty[\right)$. It is therefore natural to introduce
the following def\/inition  inspired from the notion of holomorphic family of
operators used by Kontsevich and Vishik in~\cite{KV}.

\begin{definition} \label{defn:holfamily}   We  call  a family
\begin{gather*}
A(z)\simeq \sum_{j=0}^\infty b_j(z) \vert D\vert^{\alpha(z)-j}  \in {\cal E}^{\alpha(z)}({\cal A}, D)
\end{gather*}  of
linear  operators acting on ${\cal H}^\infty$  parametrised by $z\in \C$,
 {\it holomorphic}, if
\begin{enumerate}
\itemsep=0pt
\item[1)]    $\alpha(z)$  is a  complex valued  holomorphic function,
\item[2)]    $b_j(z)$  is a ${\cal B}$-valued holomorphic function for any  non-negative integer $j$,
\item[3)] the family $A(z)\vert D\vert^{-\alpha(z)}$ is holomorphic in the Banach space ${\cal B}({\cal H})$ and
\[
\partial^{(k)}\big(A(z)\vert D\vert^{-\alpha(z)}\big)\simeq  \sum_{j=0}^\infty b_j^{(k)}(z) \vert D\vert^{-j}  \in {\cal E}^{0}({\cal A}, D).
\]
\end{enumerate}
\end{definition}

\begin{example}\label{ex:1} For any complex number $a$ and any operator  $A$ in ${\cal E}^a({\cal A},D)$, the family $A(z)= A\, \vert D\vert^{-z}$ is   a holomorphic family of order
  $\alpha(z)=a-z$ acting on ${\cal H}^\infty$.
\end{example}

The following lemma extends this to  families $A(z)=  \vert D\vert^{-z}\, A$.
\begin{lemma} For any complex number $a$ and any operator $ A\simeq \sum\limits_{j=0}^\infty
b_j \vert D\vert^{a-j}\in {\cal E}^a({\cal A}, D)$,  the product
$\vert D\vert^{-z}  A$ is a holomorphic family of order $a-z$ acting on ${\cal H}^\infty$, and the commutator $\left[\vert D\vert^{-z}, A\right]$ is a holomorphic family  of order $a-z-1$ acting on ${\cal H}^\infty$.
\end{lemma}

\begin{proof} Applying~(\ref{eq:commDalpha}), with
$\alpha$ replaced with $-z$,  to the operators $b_j  \vert D\vert^{a-j}$ for any non-negative integer $j$,  yields
\begin{gather*}
 \vert D\vert^{-z}  A \simeq  \sum_{k=0}^\infty c_{-z,k}
  \sum_{j=0}^\infty \delta^k(b_j)   \vert D\vert^{a-j-k-z}
 \simeq  \sum_{i=0}^\infty
 \left( \sum_{j=0}^i  c_{-z,i-j} \delta^{i-j}(b_j)\right)   \vert D\vert^{a-z-i}.\tag*{\qed}
\end{gather*}\renewcommand{\qed}{}
\end{proof}

We shall also need the following technical result.

\begin{proposition}\label{prop:Aprime}  Given a  holomorphic family $A(z)\simeq\sum\limits_{j=0}^\infty
b_j(z) \vert D\vert^{\alpha(z)-j}$  acting on ${\cal H}^\infty$, of non
  constant affine  order
  $\alpha(z)=a-q  z$, the higher   derivatives
 $\partial_z^{(k)}\left(  A(z) \vert D\vert^{-\alpha(z)}\right)$  lie
 in ${\cal E}^{0}({\cal A}, D)$, whereas   the higher   derivatives
 $\partial_z^{(k)} A(z)$  lie
 in ${\cal E}^{\alpha(z), k}({\cal A}, D)$.
\end{proposition}

\begin{proof} By assumption, the family  $A(z)  \vert D\vert^{-\alpha(z)}$ is
holomorphic in the Banach space ${\cal B}({\cal H})$, with derivatives
\[
\partial_z^k\big(A(z)\, \vert D\vert^{-\alpha(z)}\big)\simeq
\sum_{j=0}^\infty \partial_z^k b_j(z)   \vert D\vert^{-j}
\]
in ${\cal E}^0({\cal A}, D)$. Since  $\vert D\vert^{\alpha(z)}$ is a holomorphic
   family with derivatives
   \[
   \partial_z^l\big(\vert D\vert^{\alpha(z)}\big)=
   (-q)^l \log^l \vert D\vert  \vert D\vert^{\alpha(z)},
   \] that lie
in ${\cal E}^{\alpha(z), l}({\cal A }, D)$, the Leibniz rule yields the
following identity of (unbounded) operators:
\begin{gather*}
A^{(k)}(z)= \sum_{l=0}^k {k\choose l} (-q)^l
\left(\sum_{j=0}^\infty \partial_z^{k-l}b_j(z)   \vert D\vert^{\alpha(z)-j}\right)
\log^l \vert D\vert\in {\cal E}^{\alpha(z),k}\left({\cal A}, D\right).\tag*{\qed}
\end{gather*}
\renewcommand{\qed}{}
 \end{proof}

{\sloppy
\begin{example}\label{ex:differencelog}  For any complex number $a$ and an  operator $A\in {\cal E}^a({\cal A}, D)$,  the bracket $[A, \vert D\vert^{-z}]$ and
 the dif\/ference
 $A  \left((\vert D\vert +P)^{-z}-\vert D\vert^{-z }\right)$
  def\/ine holomorphic families acting on~${\cal H}^\infty$.
  By the results of the previous section, their derivatives at zero
   $-[A,\log \vert D\vert]$ and
   $A  \left(\log \vert D\vert -\log \left(\vert D\vert +P\right)\right)$, which are
    expected to be in $ {\cal E}^{a,1}({\cal A}, D)$ actually lie in
    $ {\cal E}^{a}({\cal A}, D)$.
\end{example}

}

\section{Traces of holomorphic families}\label{section5}

Let $\left({\cal A}, {\cal H}, D\right)$ be a regular spectral triple which
satisf\/ies assumption (H)  of the previous section. We  make two  further
assumptions\begin{itemize}
\item\itemsep=0pt
{\bf Assumption (Dk).} {\it  The poles of the meromorphic map $z\mapsto  {\rm Tr}\left(b(z) \vert
    D\vert^{-z}\right)\vert^{\rm mer} $ introduced in~\eqref{eq:TRz}
have  multiplicity $\leq  k+1$.}
\item {\bf Assumption (T).} {\it  The dimension spectrum ${\rm Sd}$ is
    invariant under translation by negative integers:}
    \[
    {\rm Sd}-\N \subset {\rm Sd}.
\]
\end{itemize}

\begin{remark} Under assumption (Dk) we have
\[
{\rm Res}^j_0 \left(  {\rm Tr}\left(b(z) \vert
    D\vert^{-z}\right)\vert^{\rm mer}\right)= 0\qquad \forall\, j> k+1,
\]
 where for
    a meromorphic function $f$ with a pole of order $N$ at zero we write in a
    neighborhood of zero:
\[
f=:\sum_{j=0}^N  \frac{{\rm Res}_0^j f}{z^j}+ O(z).
\]
\end{remark}

\begin{example} Results of Kontsevich and Vishik \cite{KV} tell us that assumptions $(Dk)$ and
  $(T)$ are satisf\/ied in the classical geometric setup   corresponding to
    Example~\ref{ex:spin}. Indeed, assumption  $(Dk)$  is satisf\/ied for $k=0$ since ${\rm
    TR}\left(b(z)  \vert D\vert^{-z}\right)$ is a meromorphic map
  with simple poles. These lie  in the set $]{-}\infty , -n]\cap
  \Z$, which is stable under translation by negative integers so that
  assumption~(T) is satisf\/ied.
\end{example}

{\it  We henceforth assume that  {\rm (H)},    {\rm (Dk)} and {\rm (T)} are satisfied.}

\begin{proposition}\label{prop:merAz} Let $A(z)\simeq
    \sum\limits_{j=0}^\infty
b_j(z) \vert D\vert^{\alpha(z)-j} $ be a holomorphic
family of operators acting on~${\cal H}^\infty$ with non-constant affine order
  $\alpha(z)$.
The map
$z\mapsto  {\rm Tr} \left( A(z)\right)$ is well-defined and holomorphic on the
half plane ${\rm Re}(\alpha(z))< -n$,  and further extends to a meromorphic
 function ${\rm Tr} \left( A(z)\right)\vert^{\rm mer}$ on
the complex plane with  poles in $\alpha^{-1}(-{\rm Sd})$ with multiplicity $\leq k+1$.
\end{proposition}

\begin{proof} Let us set $\alpha(z):= -qz+a$ for some positive real
number $q$ (a similar proof holds for negative $q$). Under assumption (H), for any non-negative integer $j$ the map
\[
 z\mapsto {\rm Tr}\big(b_j(z) \vert
  D\vert^{\alpha(z)-j}\big)= {\rm Tr}\Big(b_j(z) \vert
  D\vert^{-q \big(z-\frac{a-j}{q} \big)}\Big)
  \] is holomorphic on the
   half-plane
   \[
   {\rm Re}(z)>\frac{{\rm Re}(a)+n-j}{q} \ \Longleftrightarrow \
    {\rm Re}(\alpha(z))<-n+j.
    \]
     Under Assumption (Dk), it extends to a meromorphic map $ {\rm Tr}\big(b_j(z) \vert
  D\vert^{-q \big(z-\frac{a-j}{q}\big)}\big)\vert^{\rm mer}$
  with poles of multiplicity $\leq k+1$  in the set \[
  \frac{{\rm Sd}+a-j}{q}\subset\frac{{\rm Sd}+a}{q}= \alpha^{-1 }(-{\rm Sd})
  \] as a consequence of Assumption~(T).

 On the other hand,   given a real number $r$ and an integer $N> {\rm Re}(a)+n-rq$, under Assumption~(H),  the map  $z\mapsto {\rm Tr}(K_N(z))$ is holomorphic on the half-plane ${\rm Re}(z)>r$, where
\[
K_N(z):= A(z)-\sum_{j=0}^{N-1}b_j(z)\, \vert D\vert^{\alpha(z)-j}\simeq
\sum_{j=N}^{\infty}b_j(z)\, \vert D\vert^{\alpha(z)-j}
\] is  the remainder
operator.  Combining these observations yields a meromorphic map on the
half-plane
 ${\rm Re}(z)>r$ given by
\[
{\rm Tr} ( A(z))\vert^{\rm mer}:=  \sum_{j=0}^{N-1}
{\rm Tr} \big(b_j(z) \vert D\vert^{\alpha(z)-j}\big)\vert^{\rm mer}+ {\rm
  Tr}(K_N(z)),
   \]
   with poles in $\frac{{\rm Sd}+a}{q}= \alpha^{-1 }(-{\rm Sd})$ with multiplicity $\leq k+1$. Since $r$ can be chosen arbitrarily, this provides a meromorphic extension to the whole complex plane with poles in $\alpha^{-1 }(-{\rm Sd})$ with multiplicity $\leq k+1$. The result then
 follows. \end{proof}

Using the notation of \cite{CoMo}, we set the following def\/inition.
\begin{definition}
 For an abstract pseudodif\/ferential operator $  A$ and an integer  $j$ in $ [-1, \infty[$, let
 \begin{gather}\label{eq:taukD}
\tau^{\vert D\vert}_j(A):={\rm Res}^{j+1}_{0} {\rm
    Tr} \left( A  \vert D\vert^{-z}\right)\vert^{\rm mer}.
\end{gather}
If $j=-1$, this amounts to the f\/inite part at $z=0$:
\[
\tau^{\vert D\vert}_{-1}(A)= {\rm fp}_{z=0} {\rm
    Tr} \left( A  \vert D\vert^{-z}\right)\vert^{\rm mer}.
    \]
\end{definition}
\begin{remark}  Under assumption ${\rm (Dk)}$, we have that
$\tau^{\vert D\vert}_j(A)=0$ $ \forall\, j> k$.
\end{remark}

 The
following statement extends Proposition II.1 in~\cite{CoMo}.

\begin{theorem} \label{thm:ResAz} Let $A(z)\simeq \sum\limits_{j=0}^\infty
b_j(z) \vert D\vert^{\alpha(z)-j}$, $b_j(z)\in {\cal B}$, be a holomorphic
family of non-constant affine order
  $\alpha(z)=a-q  z$ for some positive~$q$.
\begin{enumerate}
\itemsep=0pt
\item For any integer $j\in [[-1, k]]:= \Z\cap [-1, k]$, the  $(j+1)$-th residue reads
 \begin{gather*}
{\rm Res}^{j+1}_{0} {\rm Tr}(A(z))\vert^{\rm mer}= \sum_{n=0}^{k-j}
\frac{\tau_{j+n}^{\vert D\vert}\big(\big(\partial_z^n (A(z)
  \vert D\vert^{q\, z})\big)_{\vert_{z=0}}\big)}{q^{j+n+1}  n!}.
\end{gather*}
\item In particular,
\[
{\rm fp}_{z=0}\left( {\rm Tr}(A(z))\vert^{\rm mer}\right)=
\sum_{n=0}^{k+1}
\frac{\tau_{n-1}^{\vert D\vert}\big(\big(\partial_z^n (A(z)
  \vert D\vert^{q z})\big)_{\vert_{z=0}}\big)}{q^{n}  n!},
\] where as before, ${\rm fp}_{z=0}$ stands for the finite part at $z=0$.
\end{enumerate}
\end{theorem}

\begin{proof}
1. We write $A(z)= A(z) \vert D\vert^{q\, z}\, \vert D\vert^{-q\, z}$ and implement
   a Taylor expansion on the holomorphic family $A(z) \vert
   D\vert^{-
q\, z}$, which has constant order $a$:
\[
{\rm Tr}(A(z))\vert^{\rm mer}= \sum_{n=0}^\infty \frac{z^n}{n!}  {\rm Tr}\big(\big(\left(\partial_z^n (A(z)
    \vert D\vert^{q\, z})\right)_{\vert_{z=0}} \vert D\vert^{-q
    z}\big)\vert^{\rm mer}\big).
    \]
By Proposition \ref{prop:Aprime}, the higher derivatives $\left(\partial_z^n
(A(z)  \vert D\vert^{q  z})\right)_{\vert_{z=0}}$ are operators of order $a$. We have
\begin{gather*}
{\rm Res}^{j+1}_{0}\left( {\rm Tr}(A(z))\right)\vert^{\rm mer} =  \sum_{n=0}^\infty {\rm Res}_0^{n+j+1} \frac{
\big( {\rm Tr}\big(\big(\partial_z^n  (A(z)
    \vert D\vert^{q z} )_{\vert_{z=0}}  \vert D\vert^{-q  z}\big)\big)\big)}{n!}\\
    \phantom{{\rm Res}^{j+1}_{0}\left( {\rm Tr}(A(z))\right)\vert^{\rm mer}}{}
 =  \sum_{n=0}^\infty\frac{ \tau^{\vert D\vert}_{n+j}
\big(\left(\partial_z^n \left(A(z) \vert D\vert^{q  z}\right)\right)_{\vert_{z=0}}\big)\vert^{\rm mer}}{q^{n+j+1}  n!}\\
\phantom{{\rm Res}^{j+1}_{0}\left( {\rm Tr}(A(z))\right)\vert^{\rm mer}}{}
 =  \sum_{n=0}^{k-j}\frac{ \tau^{\vert D\vert}_{n+j}\big(\left(\partial_z^n \left(A(z) \vert D\vert^{q  z}\right)\right)_{\vert_{z=0}}\big)}{q^{n+j+1}  n!},
\end{gather*}
since  for any abstract pseudodif\/ferential operator $A$, we have $\tau_{k+l}^{\vert D\vert}(A)=0$ $\forall\, l>0$.

2. The f\/inite part at $z=0$ is obtained from setting $j=-1$ in the
  previous formula.
\end{proof}
\begin{remark} These formulae are similar to formulae for the higher
  residues of traces of holomorphic families of log-polyhomogeneous
  pseudodif\/ferential operators of logarithmic type $k$,   acting on smooth sections of
a vector bundle, which were derived in \cite{L}.
\end{remark}

\section{Discrepancies}\label{section6}

{\it  As before we assume that   {\rm (H)},    {\rm (Dk)} and {\rm (T)} are satisfied.}

We want to measure the   obstructions preventing  the linear forms  $\tau_j^{\vert D\vert}$ from
having the  expected properties of a trace, which we interpret as
discrepancies of the linear forms $\tau_j^{\vert D\vert}$.

 The $\vert D\vert$-dependence of the $\vert D\vert$-weighted residue trace $\tau_j^{\vert D\vert}(A)$ of  an abstract
pseudo\-di\-f\/feren\-tial operator $A$   is one of these defects. We show how  a zero
order perturbation $\vert D\vert \to \vert D\vert + P$ of the weight $\vert
D\vert$ by some $P\in {\cal E}^0({\cal A}, D)$ af\/fects   $\vert
D\vert$-weighted residue traces. We express this variation  in terms of
expressions of the type $\tau_{j+n}^{\vert D\vert}({\cal P}_n(A,\vert
D\vert))$, where
${\cal P}_n(A,\vert D\vert)$ are abstract pseudodif\/ferential operators indexed by
positive integers $n$.

On the grounds of formulae which hold in the classical setup (see formulae (14)
and (15) in  \cite{Pa2} with $\phi(x)=x^{-z}$), similar to formulae derived in
\cite{CoMo} and \cite{Hig} (see e.g.,  the proof of Proposition~4.14),  for an operator  $P\simeq\sum\limits_{j=0}^\infty b_j \, \vert D\vert^{-j} \in {\cal E}^0({\cal A}, D)$ and a complex number
 $z$ we def\/ine $(\vert D\vert +P)^{-z}$ by
\begin{gather}\label{eq:phiDpert}
  (\vert D\vert+P)^{-z}
:\simeq \vert D\vert^{-z} \\
\qquad{}+\sum_{n=1}^\infty \sum_{\vert
  k\vert =0}^\infty (-1)^{\vert k\vert  +n} \frac{\Gamma(z+\vert k\vert
  +n)}{\Gamma(z)}
 \frac{\delta^{ k_1}(P)  \cdots  \delta^{ k_n}(P)   \vert D\vert^{-z-\vert k\vert -n}}{k!  (k_1+1)(k_1+k_2+2)\cdots (k_1+\cdots +k_{n}+n)}
,\nonumber
\end{gather}
where we have set $k!=k_1!\cdots k_n!$ and $\vert k\vert= k_1+\cdots +k_n$. Note
that the real part of the order $-z-\vert k\vert-n$ of the operators
$\delta^{ k_1}(P)  \cdots  \delta^{ k_n}(P)
    \vert D\vert^{-z-\vert k\vert -n}$ decreases with $n$ and $\vert k\vert$.
Using~(\ref{eq:commDalpha}), one can write  $(\vert D\vert+P)^{-z}$ as $\sum\limits_{j=0}^\infty \tilde b_j(z)  \vert D\vert^{-z-j}$ with $\tilde b_j(z)$, $j\in \Z_{\geq 0}$ holomorphic in ${\cal B}$,  so that  $D(z):=(\vert D\vert+P)^{-z}$  def\/ines a holomorphic family of abstract pseudodif\/ferential operators of order $-z$. It follows from Proposition~\ref{prop:merAz} that for any abstract pseudodif\/ferential operator $A$ acting on ${\cal H}^\infty$, the map $z\mapsto {\rm Tr} \left(A (\vert D\vert+P)^{-z}\right)$ is holomorphic on some half-plane and extends to a~meromorphic function ${\rm Tr} \left(A (\vert D\vert+P)^{-z}\right)\vert^{\rm mer}$ with  poles of order $\leq k+1$.
Let us set: \begin{gather*}
\tau^{\vert D\vert+P}_j(A):={\rm Res}^{j+1}_{0} {\rm
    Tr} \left( A  (\vert D\vert+P)^{-z}\right)\vert^{\rm mer}.
\end{gather*}

\begin{proposition}  Given  $P\simeq\sum\limits_{j=0}^\infty b_j   \vert D\vert^{-j}$ in $ {\cal E}^0({\cal A}, D)$, for any non-negative  integer $k$ and any
  $j\in[[-1,k-1]]$ \begin{gather}\label{eq:taujD1D2}\tau_j^{\vert D\vert +P}(A)-\tau_j^{\vert D\vert}(A)= \sum_{n=1}^{k-j}
\frac{\tau_{j+n}^{\vert D\vert}\left( A \left(\partial_z^n((\vert D\vert +P)^{-z} \vert
    D\vert^z)_{\vert_{z=0}}\right)\right)}{n!}
\end{gather}
holds for any  operator $A$ in $ {\cal E}^a({\cal A}, D)$ of any complex order $a$.
Moreover,
\[
\tau_k(A):= \tau_k^{\vert D\vert+P}(A)= \tau_k^{\vert D\vert}(A)
 \]
 is independent of the perturbation $P$ and when $k=0$ the
above
formula reads
\[
 \tau_{-1}^{\vert D\vert +P}(A)-\tau_{-1}^{\vert D\vert}(A)= \tau_{0}\left( A  \left(\log
   ( \vert D\vert+P) -\log \vert D\vert\right) \right).
   \]
\end{proposition}

\begin{proof} We  apply  Theorem \ref{thm:ResAz} to the holomorphic family
\[
A(z)= A  ((\vert D\vert+P)^{-z} -\vert D\vert^{-z})=
 A  ((\vert D\vert+P)^{-z} \vert D\vert^{z}-1)\vert D\vert^{-z},
 \] of order
 $\alpha(z)= a-z$.

 By Proposition \ref{prop:Aprime}, the  higher derivatives
at zero
\[
\left(\partial_z^n\left(A(z) \vert
    D\vert^{z}\right)\right)_{\vert_{z=0}}=  A  \partial_z^n\left((\vert
D\vert+P)^{-z} \vert D\vert^{z}\right)_{\vert_{z=0}}
\]
have order zero,  and by Theorem~\ref{thm:ResAz}, we have
\begin{gather*}
\tau_j^{\vert D\vert +P}(A)-\tau_j^{\vert D\vert}(A) =  {\rm Res}_{0}\left(z^j\, {\rm Tr}\left(
 A((\vert D\vert+P)^{-z} -\vert D\vert^{-z})\right) \vert^{\rm mer}\right)\\
 \phantom{\tau_j^{\vert D\vert +P}(A)-\tau_j^{\vert D\vert}(A)}{}
= {\rm Res}_0^{j+1}  \left( {\rm Tr}(A((\vert D\vert+P)^{-z}-\vert
D\vert^{-z}))\vert^{\rm mer}\right)
\\
 \phantom{\tau_j^{\vert D\vert +P}(A)-\tau_j^{\vert D\vert}(A)}{}
 =   \sum_{n=  0}^{k-j} \frac{\tau_{j+n}^{\vert D\vert}\left( A \left(\partial_z^n((\vert D\vert+P)^{-z}  \vert
    D\vert^z-1)_{\vert_{z=0}}\right)\right)}{n!}\\
 \phantom{\tau_j^{\vert D\vert +P}(A)-\tau_j^{\vert D\vert}(A)}{}
 =   \sum_{n=1 }^{k-j} \frac{\tau_{j+n}^{\vert D\vert}\left( A \left(\partial_z^n(\vert D\vert+P)^{-z}  \vert
    D\vert^z)_{\vert_{z=0}}\right)\right)}{n!}\qquad {\rm if}\quad j<k,
\end{gather*}
since the $n=0$   term  vanishes. If $j=k$, we have
$\tau_j^{\vert D\vert +P}(A)-\tau_j^{\vert D\vert}(A)=0$, which shows the inva\-riance of
$\tau_k^{\vert D\vert}$ under a perturbation $P$. By Proposition
 \ref{prop:differencelog},  the operator  $A  \left(\log
    (\vert D\vert+P) -\log \vert D\vert\right)$ lies in ${\cal E}^a({\cal A }, D)$
  and when
  $k=0$, the above computation yields
  \begin{gather*}
   \tau_{-1}^{\vert D\vert +P}(A)-
\tau_{-1}^{\vert D\vert}(A)= \tau_{0}\left( A  \left( \log \vert D\vert-   \log (\vert D\vert+P)\right) \right).\tag*{\qed}
\end{gather*}
\renewcommand{\qed}{}
   \end{proof}

\begin{remark}When $k=0$, the linear form $\tau_0$ generalises Wodzicki's  noncommutative
  resi\-due~\mbox{\cite{Wo1,Wo2}}
  on classical pseudodif\/ferential operators acting on smooth sections of a
  vector bund\-le~$E$ over a closed manifold $M$   (see~\cite{Ka} for
  a review).
\end{remark}
The linear form $\tau_j^{\vert D\vert}$ does not vanish on commutators as could
be expected of a trace, leading to another discrepancy. The following result yields back  Proposition II.1 in~\cite{CoMo}.

\begin{proposition} Given two abstract pseudodifferential operators $A$, $B$ we have:
\begin{gather}\label{eq:taujDbracket}\tau_j^{D}([A,B])=\sum_{n= 1}^{k-j}
 (-1)^{n+1}\frac{\tau_{j+n}^{\vert D\vert} \left(
    A L^n(B)\right)}{n!},\qquad \forall\,
  j\in[[-1,k-1]],  \end{gather}
where we have set  $L(A):=  [\log \vert D\vert, A].$

Moreover $\tau_k([A,B])=0$ and when $k=0$ the above formula reads:
\[
\tau_{-1}^{D}([A,B ])=\tau_{0}\left( A \left([ \log
    \vert D\vert ,B] \right)\right).
    \]
\end{proposition}

\begin{proof}  Using the cyclic property of the
usual trace Tr, by analytic continuation we have:
\begin{gather*}
{\rm Tr} (B  A
  \vert D\vert^{-z} )\vert^{\rm mer} = {\rm Tr}\big(\big(\vert D\vert^{-z/2}B\big)  \big(A
  \vert D\vert^{-z/2}\big)\big)\vert^{\rm mer}\\
  \phantom{{\rm Tr} (B  A   \vert D\vert^{-z} )\vert^{\rm mer} }{}
 = {\rm Tr}\big( \big(A
  \vert D\vert^{-z/2}\big) \big( \vert D\vert^{-z/2}B\big)\big)\vert^{\rm mer}
 = {\rm Tr} ( A  \vert D\vert^{-z} B )\vert^{\rm mer}.
\end{gather*} Hence,
\begin{gather*}
{\rm Tr} ( A  B\vert D\vert^{-z}-B  A
  \vert D\vert^{-z} )\vert^{\rm mer} =  {\rm Tr} ( A  B \vert
  D\vert^{-z}-A \vert D\vert^{-z}B )\vert^{\rm mer}
 =  {\rm Tr} ( A [B ,\vert
  D\vert^{-z} ] )\vert^{\rm mer} .
\end{gather*}
We can apply Theorem \ref{thm:ResAz} to the holomorphic family
\[
A(z)= A\, [B, \vert
  D\vert^{-z}].
  \] Indeed,
\[
A(z) \vert
    D\vert^z= A  B -A \vert D\vert^{-z}B \vert
    D\vert^z= A   (B- \sigma(z)(B)),
    \]
where we have set $\sigma(z)(A)= \vert D\vert^{-z}  A \vert
    D\vert^z$. Since
    \[
    \vert D\vert^{-z} A\simeq  A  \vert D\vert^{-z}+
     \sum_{k=1}^\infty c_{-z,k}
\delta^k(A)   \vert D\vert^{-z-k},
\] it follows that
\[
\sigma(z)(A)\simeq A+ \sum_{k=1}^\infty c_{-z,k}
\delta^k(A)   \vert D\vert^{ -k},
\] where as before
 $c_{ \alpha, k}= \frac{\alpha(\alpha-1)\cdots (\alpha-k+1)}{k!}$ for any positive
 integer $k$. By
 Proposition \ref{prop:Aprime}, the  higher derivatives at zero $\partial_z^n (A(z) \vert
    D\vert^z)_{\vert_{z=0}}$  have order $a+b$, where $a$ is the order of $A$
    and $b$ that of $B$.

Since  $\partial^n_z \sigma(z)(A)_{\vert_{z=0}}=
(-1)^n\, L(A)$,   Theorem \ref{thm:ResAz}
    yields
\begin{gather*}
 \tau_j^{\vert D\vert}([A,B]) =  \sum_{n= 0}^{k-j} \frac{\tau_{j+n}^{\vert D\vert} \left(\partial_z^n (A(z)
  \vert D\vert^z)\right)_{\vert_{z=0}}}{n!}
 =  -\sum_{n= 0}^{k-j} \frac{\tau_{j+n}^{\vert D\vert} \left(A\, \left(\partial_z^n \sigma(z)
  (B)\right) \right)_{\vert_{z=0}}}{n!}\\
 \phantom{\tau_j^{\vert D\vert}([A,B])}{}
 =  \sum_{n= 0}^{k-j} (-1)^{n+1}\frac{\tau_{j+n}^{\vert D\vert} \left(
    A L^n(B)\right)}{n!}
 =  \sum_{n=1}^{k-j} (-1)^{n+1}\frac{\tau_{j+n}^{\vert D\vert} \left(
    A L^n(B)\right)}{n!}\qquad{\rm if}\quad j<k
    \end{gather*}
since the $n=0$ term vanishes.

  This computation also shows that $ \tau_k([A,B]) =\tau_k^{\vert D\vert} ([A,B])=0$.  When
  $k=0$ it  shows that:
\begin{gather*}
\tau_{-1}^{\vert D\vert}([A,B ])=\tau_{0}\left( A \left([ \log
    \vert D\vert ,B] \right)\right).\tag*{\qed}
    \end{gather*}\renewcommand{\qed}{}
\end{proof}

These discrepancy formulae are similar to the anomaly formulae for weighted
traces of classical  pseudodif\/ferential operators, see \cite{MN,CDMP,PS}  and
 \cite{Pa1}  for an extension  to  log-polyhomogeneous operators.

\section{ An analog of Kontsevich and Vishik's canonical trace}\label{section7}

In \cite{KV} Kontsevich and Vishik popularised what is known as the canonical trace on
non-integer order classical pseudodif\/ferential operators acting on smooth
sections of a vector bunlde $E$ over a closed manifold $
M$. This linear form
which vanishes on commutators of non-integer order classical pseudodif\/ferential
operators, is the unique (up to a multiplicative constant) linear extension of the ordinary trace to the set of
non-integer order operators with that property. The non integrality assumption
on the order can actually be weakened to  the order not lying in the set  $[{-}n,+\infty[\cap
\Z$ which corresponds to $- {\rm Sd}$, with ${\rm Sd}$ the dimension spectrum of the spectral triple
$\left(\Ci(M), L^2(M, E), D\right).$

For a spectral triple which fulf\/ills assumptions (H), (Dk) and (T), we replace
this  assumption on the order  with the requirement that the
order lies outside the set $-{\rm Sd}$.

\begin{lemma} For an abstract pseudodifferential  operator $A$ whose order~$a$ lies outside the set $ -{\rm Sd}$  and for any
  holomorphic family $A(z)$ such that $A(0)=A$,
 the map $z\mapsto {\rm Tr} \left( A(z) \right)\vert^{\rm mer}$ is holomorphic at
  $z=0$.

 In particular,  the map $z\mapsto  {\rm
    Tr} \left( A\, \vert D\vert^{-z}\right)\vert^{\rm mer} $ is holomorphic at
  zero so that
  \[
  \tau_j^{\vert D\vert}(A)=0\qquad \forall\, j\in \Z_+\qquad {\rm and}\qquad \tau_{-1}^{\vert D\vert} (P)=\lim_{z\to 0} ({\rm
    Tr}  ( A \vert D\vert^{-z} )\vert^{\rm mer} ).
    \]
\end{lemma}
\begin{proof} In view of the asymptotic expansion $A(z)\simeq \sum\limits_{j=0}^\infty
b_j (z)\vert D\vert^{\alpha(z)-j},$
where   $\alpha(z)$ is the  order of the holomorphic family $A(z)$,
it suf\/f\/ices to show the result for each holomorphic family $A_j(z)= b_j(z)
\vert D\vert^{\alpha(z)-j}$, with $j$ varying in $\Z_{\geq 0}$.
Such a family  has order $\alpha_j(z)= \alpha(z)-j$ so that the  poles of the meromorphic map $z\mapsto {\rm
    Tr} \left( A_j(z) \right)\vert^{\rm mer}$  lie in
  $\alpha^{-1}(-{\rm Sd}+\Z_{\geq 0})\subset \alpha^{-1}(-{\rm Sd})$, where we have used
  assumption (T) on the spectrum.
 Since    the order of $A$ does not lie in $-{\rm Sd}$, $0$ does not lie in this set of
poles, which leads to the two identities in the statement.\end{proof}

\begin{proposition}\label{proposition6}   For  an abstract pseudodifferential  operator $A$ whose order
  $a$ does not lie in~$-{\rm Sd}$,
\[
\tau_{-1}(A):=\tau_{-1}^{\vert D\vert} (A)
\] is insensitive to lower
order perturbations $\vert D\vert\to \vert D\vert +P$ of the weight
 $\vert D\vert$,  with $P\simeq \sum\limits_{j=0}^\infty b_j   \vert D\vert^{-j}
 $ in  ${\cal E}^0({\cal A}, D)$.
\end{proposition}
\begin{proof}
By (\ref{eq:taujD1D2}) and for $j$ in $[[-1, k-1]]$
\begin{gather*}
 \tau_{-1}^{\vert D\vert +P}(A)-\tau_{-1}^{\vert D\vert}(A) =  \sum_{n=1}^{k+1}
\frac{\tau_{-1+n}^{\vert D\vert}\left( A \left(\partial_z^n((\vert D\vert +P)^{-z} \vert
    D\vert^z)_{\vert_{z=0}}\right)\right)}{n!}\\
\phantom{\tau_{-1}^{\vert D\vert +P}(A)-\tau_{-1}^{\vert D\vert}(A)}{}
 =  \sum_{m= 0}^k
\frac{\tau_{m}^{\vert D\vert}\left( A \left(\partial_z^{m+1}((\vert D\vert+P)^{-z} \vert
    D\vert^z)_{\vert_{z=0}}\right)\right)}{(m+1)!}.
\end{gather*}
By Proposition \ref{prop:Aprime}, the operators  $ A \left(\partial_z^{m+1}((\vert D\vert+P)^{-z} \vert
    D\vert^z)\vert_{z=0}\right)$, $m\in \Z_{\geq 0}$, all  have orders $a\notin -{\rm Sd}$, as a consequence
  of which
   all the terms $\tau_{m}^{\vert D\vert}\left( A \left(\partial_z^{m+1}((\vert D\vert+P)^{-z} \vert
    D\vert^z)_{\vert_{z=0}}\right)\right)$  in the sum
vanish. Thus,
\begin{gather*}
\tau_{-1}^{\vert D\vert +P}(A)=\tau_{-1}^{\vert D\vert}(A).\tag*{\qed}
\end{gather*}\renewcommand{\qed}{}
\end{proof}

The linear form $\tau_{-1}$  provides a substitute on
abstract pseudodif\/ferential operators with order outside the discrete set $-{\rm Sd}$, for the canonical trace TR
def\/ined on classical pseudodif\/ferential
operators, whose order lies outside  the set $[{-}n,+\infty[\cap \Z$.
\begin{proposition}\label{proposition7}
The linear\footnote{Meaning by this that it acts on  linear
  combinations of operators whose order does not lie in $-{\rm Sd}$.}  map
$\tau_{-1}$  vanishes on  commutators, whose order does not  lie in~$-{\rm Sd}$:
\[
\left({\rm ord}(A)+{\rm ord}(B)\notin -{\rm Sd}\right) \ \Longrightarrow \
\tau_{-1}\left([A,B]\right)=0.
\]
\end{proposition}

\begin{proof}
By (\ref{eq:taujDbracket}) we have
\begin{gather*}
\tau_{-1}([A,B]) = \sum_{n> 0} (-1)^{n+1}\frac{\tau_{n-1}^{\vert D\vert} \left(
    A L^n(B)\right)}{n!}
 =  \sum_{m\geq  0} (-1)^{m}\frac{\tau_{m}^{\vert D\vert} \left(
    A L^{m+1}(B)\right)}{(m+1)!}.
\end{gather*}
The orders of the  operators
  $  A  L^{m+1}(B)$, which lie  in ${\rm ord}(A)+ {\rm ord}(B)-\Z_{\geq 0}$  do not lie
  in $-{\rm Sd}$ since   ${\rm ord}(A)+ {\rm ord}(B)$ does not and since by assumption (T) we have the inclusion: $ -{\rm Sd}\subset
-{\rm Sd} -\Z_{\geq 0}$. Hence,
 the terms $\tau_{m}^{\vert D\vert} \left(
    A\,L^{m+1}(B)\right)$, $m\in \Z_{\geq 0}$ in the sum vanish and
$\tau_{-1}([A,B])=0$. \end{proof}

 Thus, the linear map $\tau_{-1}$ def\/ines a ``trace''  on the set of
operators, whose order  lies outside the set $-{\rm Sd}$  in the same way that
Kontsevich and Vishik's canonical trace TR
def\/ines a~``trace'' on non-integer order classical pseudodif\/ferential
operators.

\subsection*{Acknowledgements}  I thank the Australian National University in
Canberra for its hospitality during a ten days stay, when this work was
initiated and  Bai Ling Wang for inviting me there, as well as for stimulating
workshops he organised and the lively discussions that followed. I am grateful
to Alan Carey for helping organise my stay and I thank him and Adam Rennie
most warmly  for
various  informal discussions we had during my stay, which triggered this work. I would also like to acknowledge the referees'
substantial help in improving 
of this paper.

\pdfbookmark[1]{References}{ref}
\LastPageEnding

\end{document}